\documentclass[letterpaper, draftclsnofoot, onecolumn]{ieeeconf}

\IEEEoverridecommandlockouts                              

\usepackage[pdftex]{graphicx}
\usepackage{amsmath}
\usepackage{amssymb}  
\usepackage{cite}
\usepackage{array}
\usepackage{color}
\usepackage{url}
\usepackage{subfig}
\usepackage{stackrel}
\usepackage{algorithm,tabularx}
\usepackage[noend]{algpseudocode}
\usepackage{mathtools}

\makeatletter
\newcommand{\multiline}[1]{%
  \begin{tabularx}{\dimexpr\linewidth-\ALG@thistlm}[t]{@{}X@{}}
    #1
  \end{tabularx}
}
\makeatother

\newcommand{\real}{{\mathbb{R}}}
\newcommand{\subj}{\text{subj. to}}

\newcommand\oprocendsymbol{\hbox{$\square$}}
\newcommand\oprocend{\relax\ifmmode\else\unskip\hfill\fi\oprocendsymbol}

\newcommand{\perfHinf}{\nu}

\newcommand{\xlin}{\Delta x}
\newcommand{\ulin}{\Delta u}

\newcommand{\xdir}{z}
\newcommand{\udir}{v}

\newcommand{\xproj}{\varphi}
\newcommand{\uproj}{\gamma}

\newcommand{\finalTime}{T}

\newcommand{\Tzm}{\in \mathbb{T}_{[0,\finalTime]}}
\newcommand{\Tzmo}{\in \mathbb{T}_{[0,\finalTime-1]}}
\newcommand{\Tom}{\in \mathbb{T}_{[1,\finalTime]}}

\newcommand{\indexTv}{s}
\newcommand{\numTv}{S}
\newcommand{\indexVertex}{p}
\newcommand{\numVertex}{P}

\newcommand{\sparseMatrixVertex}{Adj_\indexVertex}
\newcommand{\sparseMatrixVertexElement}[1]{adj_{\indexVertex#1}}

\newcommand{\at}[1]{a_{t#1}^k}  

\newcommand{\amin}[1]{a_{#1}^{min}}
\newcommand{\amax}[1]{a_{#1}^{max}}

\newcommand{\Av}{\tilde{A}}
\newcommand{\av}{\tilde{a}}
\newcommand{\Avertex}{\Av_\indexVertex} 
\newcommand{\AvertexElement}[1]{\av_{\indexVertex#1}}

\newcommand{\Kv}{\tilde{K}}
\newcommand{\Kvertex}{\Kv_\indexVertex} 

\newcommand{\coeffvertex}{\theta_{\indexVertex,t}^k}

\newcommand{\Amin}{\underbar{$A$}}
\newcommand{\Amax}{\bar{A}}
\newcommand{\AminElement}{\underbar{$a$}}
\newcommand{\AmaxElement}{\bar{a}}

\newcommand{\Kmin}{\underbar{$K$}}
\newcommand{\Kmax}{\bar{K}}
\newcommand{\KminElement}{\underbar{$k$}}
\newcommand{\KmaxElement}{\bar{k}}
\newcommand{\GG}{{\mathcal{G}}}
\newcommand{\NN}{{\mathcal{N}}}
\newcommand{\EE}{{\mathcal{E}}}

\newtheorem{theorem}{Theorem}
\newtheorem{remark}{Remark}
\newtheorem{assumption}{Assumption}
\newtheorem{proposition}{Proposition}

\title{\LARGE \bf
A Sparse Polytopic LPV Controller\\ for Fully-Distributed Nonlinear Optimal Control
}

\author{Sara Spedicato, Sarnavi Mahesh and Giuseppe Notarstefano
\thanks{Sarnavi Mahesh is with the
Department of Engineering, Universit\`a del Salento, Via per Monteroni, Lecce, Italy,
        {\tt\small sarnavi.mahesh@unisalento.it}}%
\thanks{Sara Spedicato and Giuseppe Notarstefano are with the 
Department of Electrical, Electronic and Information Engineering, University of Bologna, 
Viale del Risorgimento 2, Bologna, Italy,
        {\tt\small name.lastname@unibo.it}}%
}

\begin{document}

\maketitle
\thispagestyle{empty}
\pagestyle{empty}


\textbf{
© 2019 IEEE.  Personal use of this material is permitted.  Permission from IEEE must be obtained for all other uses, in any current or future media, including reprinting/republishing this material for advertising or promotional purposes, creating new collective works, for resale or redistribution to servers or lists, or reuse of any copyrighted component of this work in other works.
}

\vspace{0,5cm}

\begin{abstract}
  In this paper we deal with distributed optimal control for nonlinear dynamical
  systems over graph, that is large-scale systems in which the dynamics of each
  subsystem depends on neighboring states only. Starting from a previous work in
  which we designed a partially distributed solution based on a cloud, here we
  propose a fully-distributed algorithm. The key novelty of the approach in this
  paper is the design of a sparse controller to stabilize trajectories of the
  nonlinear system at each iteration of the distributed algorithm. The proposed
  controller is based on the design of a stabilizing controller for polytopic
  Linear Parameter Varying (LPV) systems satisfying nonconvex sparsity
  constraints.  Thanks to a suitable choice of vertex matrices and to an
  iterative procedure using convex approximations of the nonconvex matrix
  problem, we are able to design a controller in which each agent can locally
  compute the feedback gains at each iteration by simply combining coefficients
  of some vertex matrices that can be pre-computed offline. We show the
  effectiveness of the strategy on simulations performed on a multi-agent
  formation control problem.
\end{abstract}

\section{Introduction}
\emph{Nonlinear} optimal control of network systems is a challenging problem
with applications in several control areas as cooperative robotics, smart grids
or spatially distributed control systems.
The large-scale nature and nonconvexity of the optimization problem
are the main challenges that need to be taken into account in addressing the
solution of these optimal control problems in a distributed way.

Distributed optimal control over networks has been mainly investigated for
linear (time-invariant) systems, \cite{doan2011iterative,
  giselsson2013accelerated, conte2016distributed, gross2016cooperative}, so that
the resulting optimization problem is convex.  While the approaches developed
for convex problems are fully distributed, in the few methods proposed for
nonlinear, nonconvex problems, \cite{necoara2009distributed, dinh2013dual}, only
part of the computation is performed locally by the agents.
In our previous work \cite{spedicato2018necsys} we have proposed a
cloud-assisted distributed algorithm to solve (nonconvex) optimal control problems for
nonlinear dynamics over graph \cite{ferrari2006analysis}, i.e., large-scale
systems characterized by a dynamic coupling (modeled by a graph) among
subsystems. 
The algorithm proposed in \cite{spedicato2018necsys} combines distributed
computation steps with centralized steps performed by a cloud. 
A key distinctive feature of the optimal control strategy in
\cite{spedicato2018necsys} and its distributed version proposed in this paper,
is that at each iteration agents compute a trajectory of the dynamical system,
i.e., a state-input curve satisfying the dynamics. This feature is extremely
important in realtime control schemes (as, e.g., Model Predictive Control
ones) since it allows agents to stop the algorithm at any iteration and yet have
a (suboptimal) system trajectory.
This property is guaranteed through a nonlinear feedback controller acting as a
\emph{projection operator},~\cite{hauser2002projection}, that at each iteration
of the algorithm projects infeasible state-input curves into trajectories of the
system (satisfying the dynamics). 
Here, we are interested in designing a \emph{distributed} projection
operator,\cite{spedicato2018necsys}, as a sparse feedback matrix that
exponentially stabilizes the linearization of the system at a given trajectory.
The design of sparse controllers for network systems, which allows for
distributed control laws, has been investigated in the literature mainly for
linear systems. 
The works in \cite{massioni2009distributed, lin2011augmented, wu2016input, maartensson2012scalable,polyak2013lmi,jain2017online}
address the design of a sparse static feedback for linear time-invariant systems.
Dynamic controllers are instead considered in \cite{farhood2015distributed, farhood2015distributedLPV,
  viccione2009lpv, eichler2014robust}.
Among these works,\cite{farhood2015distributed} deals with time-varying systems
while \cite{viccione2009lpv, eichler2014robust,farhood2015distributedLPV}
consider LPV systems.
%

The main contributions of this paper are as follows. 
We propose a variation of the strategy introduced in \cite{spedicato2018necsys}
that allows agents to solve nonlinear optimal control problems for dynamics over
graph in a fully-distributed way, i.e., without requiring the presence of a
cloud.
In the strategy in \cite{spedicato2018necsys} the cloud computes, at
each iteration of the algorithm, a new stabilizing feedback matrix for the
current trajectory.
In this paper we remove this centralized step, so that both the implementation of
the control law and its design are purely distributed at each iteration of the
optimal control algorithm. 
The main idea is to split the design of the time-varying (and iteration
dependent) feedback controller in two parts: a computationally expensive step
performed offline before the optimal control algorithm starts, and a
computationally inexpensive one that agents perform at each iteration in a
completely distributed way.
Specifically, we are able to express each sparse time-varying, and iteration
dependent, feedback matrix as a convex combination of given vertex matrices of a
polytopic LPV system.
By simultaneously imposing suitable sparsity structures on the vertex matrices
and stability conditions for polytopic systems, we are able to obtain time and
iteration dependent feedback matrices: (i) exponentially stabilizing the
trajectories of the nonlinear system, and (ii) whose (time and iteration
dependent) convex combinators can be computed locally by agents.  In order to
obtain vertex matrices satisfying the required stability and sparsity
conditions, we propose an iterative algorithm, based on the convexification of
the nonconvex sparsity constraint, inspired by the one proposed in
\cite{fardad2014design} for linear time-invariant systems.
We extend the approach in \cite{fardad2014design} to polytopic LPV systems and
properly tailor it to distributed controllers for network systems. Also, as
opposed to \cite{viccione2009lpv, eichler2014robust,farhood2015distributedLPV},
our controller is static and can be implemented via an online step with
precomputed vertex matrices.

The paper is organized as follows. In Section \ref{sec:problem-cloud-assisted}
we describe the nonlinear optimal control problem addressed in the paper and
recall our cloud-assisted distributed algorithm proposed in
\cite{spedicato2018necsys}. In Section~\ref{sec:sparse_controller} we present
our strategy for the computation of the distributed projection operator over
graph by means of a sparse polytopic LPV controller, while in Section
\ref{sec:fully-distributed} we show how the strategy is used to obtain a fully
distributed optimal control algorithm. Finally, in
Section~\ref{sec:numerical_computations} we provide numerical computations.

\section{Problem Set-up and\\ Cloud-Assisted Algorithm}
\label{sec:problem-cloud-assisted}
\subsection{Distributed Nonlinear Optimal Control over Graph}
Nonlinear dynamics over graph consist of $N$ subsystems whose local dynamics
depend only on neighboring subsystems. The neighboring structure is modeled by a
fixed, connected and undirected graph $\GG = \{\{ 1, \dots, N\},\EE \}$, with
$\EE$ being the set of edges. We let
$\NN_i := \{ j \in \{ 1, \dots, N\} \vert (i,j) \in \EE \}$ be the set of
neighbors of node $i$, and let $Adj \in \real^{N \times N}$ denote the adjacency
matrix associated to $\GG$.  We will consider the evolution of the dynamics over a
time horizon $\finalTime$. For $0\leq t_1\leq t_2\leq \finalTime$ we define
$\mathbb{T}_{[t_1,t_2]} := \{t_1,\dots,t_2\}$.
For such a dynamical system we want to solve the nonlinear optimal control problem
\begin{align}
\min_{{\substack{x_{i,1}, \dots, x_{i,\finalTime}\\ u_{i,0}, \dots, u_{i,\finalTime-1}\\i\in \{1,\dots,N\}}}}
 \;\;& \sum_{i = 1}^{N} 
\Big( \sum_{t=0}^{\finalTime-1} \Big(\ell_i(x_{i,t},u_{i,t}) \Big) + m_i(x_{i,\finalTime}) \Big), \label{eq:cost}\\
\subj \;\; & x_{i,t+1} = f_i(x_{\mathcal{N}_i,t}, u_{i,t}), \quad t \Tzmo, \label{eq:dynamics}\\
\;\; & i \in \{1,...,N\}, \nonumber
\end{align}
where $x_{i,t} \in \real, u_{i,t} \in
\real$ are, respectively, the state and input of agent $i$ at time
$t$, $x_{i,0}$ is a (given) initial condition, $x_{\mathcal{N}_i,t} \in
\real^{|\NN_i |}$, where $|\NN_i|$ is the cardinality of
$\NN_i$, is a vector with components $x_{j,t}$, $j \in \mathcal{N}_i$, $\ell_i :
\real \times \real \rightarrow \real$, $m_i : \real \rightarrow
\real$ are cost functions, and $f_i : \real^{|\NN_i |} \times \real \rightarrow
\real$ is the local state function of agent
$i$.  
The functions $\ell_i(\cdot,\cdot)$, $m_i(\cdot)$ and
$f_i(\cdot,\cdot)$, for all $i \in \{
1,\dots,N\}$, are continuously differentiable functions.  
In order to simplify the presentation of the LPV-based controller technique in
Section~\ref{sec:sparse_controller}, we suppose that the gradient of the
function $f_i(\cdot,\cdot)$ with respect to the input
$u_{i,t}$ does not depend on time $t$, for all
$i$.

A \emph{trajectory} of \eqref{eq:dynamics} consists of states and inputs,
respectivelly $x_{i,t}, \; t \Tzm, \; i \in \{ 1, \dots, N\}$ and
$u_{i,t}, \; t \Tzmo, \; i \in \{ 1, \dots, N\},$ that satisfy the dynamics
\eqref{eq:dynamics}.  Since the problem \eqref{eq:cost}-\eqref{eq:dynamics} is
nonconvex we seek for trajectories, namely $x_{i,t}^*, \; t \Tzm,$
$u_{i,t}^*, \; t \Tzmo$, $i \in \{ 1, \dots, N\}$, satisfying the dynamics,
together with Lagrange multipliers $p_{i,t}^*, \; t \Tom$,
$i \in \{ 1, \dots, N\}$, all satisfying the first-order necessary conditions
for optimality.
In our distributed scenario, agent $i$ only knows its functions
$\ell_i(\cdot,\cdot), m_i(\cdot)$ and $f_i(\cdot,\cdot)$, has computation
capabilities and communicates only with its neighbors $j \in \NN_i$. We aim to
design a distributed algorithm to solve problem
\eqref{eq:cost}-\eqref{eq:dynamics}, in which each agent $i$ aims to locally
compute its own $x_{i,t}^*$, $t \Tzm$, $u_{i,t}^*$, $t \Tzmo$, $p_{i,t}^*$,
$t \Tom$.

\subsection{Cloud-assisted distributed algorithm }
\label{sec:cloud-assisted-algorithm}
In \cite{spedicato2018necsys} we have introduced a cloud-assisted distributed
algorithm to solve problem \eqref{eq:cost}-\eqref{eq:dynamics}.  In order to
understand how such an algorithm can be fully distributed, we briefly recall its
main idea and its steps.

The algorithm is a descent method in which at each iteration agents find a local
descent direction (a state-input curve) and compute a new trajectory (satisfying
the system dynamics) through a feedback controller.
The distributed controller, which we called \emph{distributed projection
  operator over graph,} is a key element of the method (inspired by the
centralized approach in~\cite{hauser2002projection}) since it guarantees to
obtain a trajectory (satisfying the dynamics) at each iteration.
Let $\alpha_{i,t} \in \real$, $t \Tzm$, $\mu_{i,t} \in \real$, $t \Tzmo$,
$i \in \{1,\dots,N\}$, be a generic state-input curve (not satisfying the
dynamics in general) which lies in a neighborhood of a trajectory
$\tilde{x}_{i,t}, \; t \Tzm,$ $\tilde{u}_{i,t}, \; t \Tzmo$,
$i \in \{ 1, \dots, N\}$, of the nonlinear system \eqref{eq:dynamics}.
A distributed projection operator over graph, mapping the curve 
$\alpha_{i,t}, \mu_{i,t}$, for all $i$ and $t$,
into a trajectory 
$x_{i,t}, u_{i,t}$, for all $i$ and $t$, of \eqref{eq:dynamics}, is defined as
the feedback system, with $x_{i,0} = \alpha_{i,0}$,
\begin{align}
\begin{split}
x_{i,t+1} &= f_i(x_{\NN_i,t}, u_{i,t}),\\
u_{i,t} & = \mu_{i,t} + \sum_{j = 1}^N k_{t(i,j)} \Big(\alpha_{j,t} -
x_{j,t}\Big), 
\label{eq:proj_operator}
\end{split}
\end{align}
for all $i \in \{1,...,N\}$ and $t \Tzmo$, where $k_{t(i,j)} \in \real$ is the
element $i,j$ of a controller matrix $K_t \in \real^{N \times N}$ with the
following features. 
First, it has a stabilizability-like property, namely it exponentially
stabilizes the trajectory
$\tilde{x}_{i,t}, \; t \Tzm, \; \tilde{u}_{i,t}, \; t \Tzmo, \; i \in \{1,\dots,
N\}$, as $\finalTime \rightarrow \infty$. Second, it satisfies the sparsity
condition $k_{t(i,j)} = 0$ if $j \notin \NN_i$, for all $i \in \{ 1,\dots,N\}$
and $t \Tzmo$.  

In order to use a compact notation, let us denote a state-input trajectory as
$(x, u)$, where $x \in \real^{N(\finalTime+1)}$ and $u \in \real^{N\finalTime}$
are respectively the stacks of $x_{i,t}$ and $u_{i,t}$, for all $i$ and
$t$. Consistently we will use the notation $(\alpha,\mu)$ for a state-input curve.
A trajectory $(x, u)$ of \eqref{eq:dynamics} can be written, by means of
\eqref{eq:proj_operator}, as a function of a curve $(\alpha,\mu)$, i.e.,
\begin{align}
x = \xproj(\alpha,\mu),\;
u = \uproj(\alpha,\mu),
\label{eq:proj_function}
\end{align}
with suitably defined functions $\xproj(\cdot,\cdot)$ and $\uproj(\cdot,\cdot)$.
By means of \eqref{eq:proj_function} and defining
\begin{align*}
g(x,u) := \sum_{i = 1}^{N} \Big( \sum_{t=0}^{\finalTime-1} \Big(\ell_i(x_{i,t},u_{i,t}) \Big) + m_i(x_{i,\finalTime}) \Big),
\end{align*}
the dynamically constrained problem \eqref{eq:cost}-\eqref{eq:dynamics} can be written as the unconstrained problem 
\begin{align}
\min_{\alpha,\mu} \;\; g(\xproj(\alpha,\mu), \uproj(\alpha,\mu)).
\label{eq:prob_unconstrained}
\end{align}
The cloud-assisted distributed algorithm in~\cite{spedicato2018necsys}, recalled
in the next table (Algorithm \ref{alg:cloud_assisted}) from the perspective of
agent $i$, is based on a steepest descent method for the unconstrained problem
\eqref{eq:prob_unconstrained} in which trajectories are obtained through the
projection operator (designed by the cloud at each iteration).

\begin{algorithm}[H]
\caption{Cloud-assisted distributed algorithm \cite{spedicato2018necsys}}
\label{alg:cloud_assisted}
\begin{algorithmic}
\Require $x_{j,t}^{0}, u_{i,t}^{0}$, for all $t$, for $j \in \NN_i$, such that $(x^{0}, u^{0})$ is a trajectory of \eqref{eq:dynamics}
\For{$k = 0, 1, 2 \ldots$}
\State $Send2Cloud(x_{i,t}^{k}, t \Tzm, u_{i,t}^{k}, t \Tzmo)$
\State $\small{ReceiveFromCloud(\!k_{t(i,j)}^k, k_{t(j,i)}^k, j\!\in\!\NN_i, t\!\Tzmo\!)}$
\State set $p_{i,\finalTime}^k = \nabla m_i(x_{i,\finalTime}^k)$
\For{$t = \finalTime-1, \dots, 0$}
\State 
\vspace{-0.5cm}
\begin{align}
\udir_{i,t}^k = - \Big( \ell_{u,i,t}^k + b_{(i,i)} \; p_{i,t+1}^k \Big)
\label{eq:descent_v}
\end{align}
\vspace{-0.4cm}
\State receive $
a_{t(j,i)}^k, b_{(j,j)}, v_{j,t}^k, \ell_{u,j,t}^k, p_{j,t+1}^k$, $j \in \NN_i \!\setminus \!\{i\}$
\State
\vspace{-0.3cm}
\begin{small}
\begin{align}
\begin{split}
\xdir_{i,t}^k &= - \sum_{j \in \NN_i}\! \Big( \! k_{t(j,i)}^{k} v_{j,t}^k \Big)\\
p_{i,t}^k &= \!\!\!\sum_{j \in \NN_i}\!\!\! \bigg(\! \Big(a_{t(j,i)}^k - b_{(j,j)} k_{t(j,i)}^k\Big) p_{j,t+1}^k \!\!-\! k_{t(j,i)}^{k} \; \ell_{u,j,t}^k\bigg) \!+\! \ell_{x,i,t}^k
\label{eq:descent_zp}
\end{split}
\end{align}
\end{small}
\EndFor
\vspace{-0.2cm}
\State $Send2Cloud(\xdir_{i,t}^k, t \Tzm, \udir_{i,t}^k, t \Tzmo)$
\State $ReceiveFromCloud(\beta^k)$
\For{$t = 0, 1, \ldots, \finalTime-1$} 
\State receive  $\xdir_{j,t}^k, x_{j,t}^{k+1}$, $j \in \NN_i \setminus \{i\},$ 
\State update curve
\begin{align}
\begin{split}
\alpha_{j,t}^{k+1} &= x_{j,t}^k + \beta^k \xdir_{j,t}^k, \quad j \in \NN_i\\
\mu_{i,t}^{k+1} &= u_{i,t}^k + \beta^k \udir_{i,t}^k
\label{eq:curve_update}
\end{split}
\end{align}
\State update trajectory
{
\begin{align}
\begin{split}
u_{i,t}^{k+1} &= \mu_{i,t}^{k+1} + \sum_{j \in \mathcal{N}_i} k_{t(i,j)}^k \Big( \alpha_{j,t}^{k+1} - x_{j,t}^{k+1} \Big) \\
x_{i,t+1}^{k+1} &= f_i(x_{\mathcal{N}_i,t}^{k+1}, u_{i,t}^{k+1})
\label{eq:closed_loop_update}
\end{split}
\end{align}
}
\EndFor
\EndFor
\end{algorithmic}
\end{algorithm}

Let, for a generic scalar function $h(\cdot,\cdot)$, $\nabla_x h(x^k,y^k)$ and
$\nabla_y h(x^k,y^k)$ respectively denote the gradients with respect to $x$ and
$y$ evaluated at $x^k, y^k$.  At each iteration $k$, agent $i$ performs the
following steps.  First, it computes its components
$\xdir_{i,t}^k, \udir_{i,t}^k,\; t \Tzmo,$ of the whole descent direction via
\eqref{eq:descent_v}-\eqref{eq:descent_zp}, where the scalars
$\ell_{x,i,t}^k, \ell_{u,i,t}^k, a_{t(i,j)}^k, b_{(i,i)}$ are defined as
$\ell_{x,i,t}^k\! := \!\nabla_{x_{i,t}} \ell_i(x_{i,t}^k, u_{i,t}^k),
\ell_{u,i,t}^k \!:= \!\nabla_{u_{i,t}} \ell_i(x_{i,t}^k, u_{i,t}^k)$,
\begin{align}
\begin{split}
a_{t(i,j)}^k &:= \nabla_{x_{j,t}} f_i(x_{\NN_i,t}^k, u_{i,t}^k),\\
b_{(i,i)} &:= \nabla_{u_{i,t}} f_i(x_{\NN_i,t}^k, u_{i,t}^k).
\label{eq:ab} 
\end{split}
\end{align}
Second, agent $i$ performs a local curve update in which it only computes
$\alpha_{j,t}^{k+1}, j \in \NN_i, \; \mu_{i,t}^{k+1}$, for all $t$, of the
overall new curve iterate $\alpha^{k+1}, \mu^{k+1}$ via \eqref{eq:curve_update},
where $\beta^k$ is a stepsize computed by the cloud.
Third, agent $i$ executes a local trajectory update via
\eqref{eq:closed_loop_update}, in which only the $i$-th states $x^{k+1}_{i,t}$
and inputs $u^{k+1}_{i,t}$, for all $t$, are computed by means of the
distributed projection operator.
The descent direction \eqref{eq:descent_v}-\eqref{eq:descent_zp} and the
trajectory update \eqref{eq:closed_loop_update} require, respectively, the
elements $k_{t(j,i)}^{k}, \; j \in \NN_i,$ and $k_{t(i,j)}^{k}, \; j \in \NN_i,$
of the matrix $K_t^k$,
%
which is computed, at each iteration $k$, by the cloud. The cloud receives the
current $x_{i,t}^k, u_{i,t}^k$, for all $t$, from all the $N$ agents and sends
back to each agent the corresponding elements
$k_{t(j,i)}^{k}, k_{t(i,j)}^{k}, \; j \in \NN_i,$ of the matrix
$K_t^k, \; t \Tzmo$.
$Send2Cloud(\cdot)$ and $ReceiveFromCloud(\cdot)$ indicate the communication
between each agent $i$ and the cloud.

\section{Sparse Polytopic LPV Controller}
\label{sec:sparse_controller}
In this section we propose a \emph{novel distributed strategy} to compute, at
each iteration $k$ of Algorithm~\ref{alg:cloud_assisted}, a feedback matrix
$K_t^k, \; t \Tzmo,$ in a neighborhood of the trajectory iterate
$x^k_{i,t}, u^k_{i,t}$, for all $i$ and $t$.  
We recall that, for each iteration $k$, the feedback matrix $K_t^k, \; t \Tzmo,$
has to:
(i) stabilize the trajectory
$x^k_{i,t}$, $u^k_{i,t}$ (for all $i$ and $t$) of \eqref{eq:dynamics} as
$\finalTime \rightarrow \infty$, and (ii) satisfy the sparsity condition
$k_{t(i,j)}^k = 0, \text{ if } j \notin \NN_i, i \in \{ 1,\dots,N\}, t \Tzmo.$

In order to satisfy (i), we use the following property. A feedback that
exponentially stabilizes the linearization of the system at a given trajectory
also locally exponentially stabilizes the trajectory of the nonlinear system.
As for (ii), let $Adj^c := \boldsymbol{1} - Adj$, with $\boldsymbol{1}$ the
matrix with all entries equal to one and $Adj$ the adjacency matrix, with
elements $adj_{(i,j)} = 1$ if $j \in \NN_i$ and $0$ otherwise.
The sparsity condition (ii) can be written compactly as
$K_t^k \circ Adj^c = 0, \quad t \Tzmo$, where $\circ$ denotes element-wise
multiplication.

We can, thus, pose the problem of designing $K_t^k, t \Tzmo,$ satisfying (i) and
(ii) as follows.
Let us consider the linearization of the nonlinear dynamics \eqref{eq:dynamics}
at the trajectory $x_{i,t}^k, u_{i,t}^k$, for all $i$ and $t$, i.e.,
\begin{align}
\begin{split}
\xlin_{t+1} &=  A_t^k \xlin_t + B \ulin_t, \quad t  \Tzmo,
\label{eq:system_z1}
\end{split}
\end{align}
where, $\xlin_t \in \real^{N}$ and $\ulin_t \in \real^{N}$ are, respectively,
state and input of the linearization system at time $t$,
$A_t^k \in \real^{N \times N}$ and $B \in \real^{N \times N}$ are, respectively,
the matrices with non zero elements
$\at{(i,j)}, \; i \in \{1,\dots, N \}, \; j \in \NN_i,$ and
$b_{(i,i)}, \; i \in \{1,\dots, N \}$, defined in \eqref{eq:ab}.  We aim to
design, at each iteration $k$ of the algorithm, a control law
\begin{align}
\ulin_t = - K_t^k \xlin_t, \quad t  \Tzmo,
\label{eq:control_law}
\end{align}
that stabilizes the $k$-th system \eqref{eq:system_z1} as $\finalTime \rightarrow \infty$
and satisfies the sparsity condition
\begin{align}
K_t^k &\circ Adj^c = 0, \quad t  \Tzmo.
\label{eq:sparsity_distributed}
\end{align}

\begin{remark}
  We consider, for simplicity a constant $B$ but the strategy we propose can be
  applied with suitable modifications to the case of a time-dependent
  matrix.~\oprocend
\end{remark}

\subsection{Main idea for the controller design}
The main idea to compute feedback matrices $K_t^k, \; t \Tzmo$ at each iteration
$k$ of the algorithm is the following.

First, in Section \ref{sec:polytopic_representation}, we show that system
\eqref{eq:system_z1} can be written as a polytopic LPV system. That is,
$A_t^k, \; t \Tzmo$, for all $k$, can be written as
\begin{align}
A_t^k = \sum_{\indexVertex = 1}^{\numVertex} \coeffvertex \Avertex,
\label{eq:At_convex}
\end{align}
where $P$ is the number of vertices of the polytope,
$\coeffvertex \in \real, \; \indexVertex = 1, \dots, \numVertex$, are suitable
vertex coefficients satisfying,
\begin{align}
\coeffvertex \geq 0, 
\quad \text{and} \quad \sum_{\indexVertex = 1}^{\numVertex}  \coeffvertex = 1,
\label{eq:coefficients}
\end{align}
and $\Avertex \in \real^{N \times N}, \; \indexVertex = 1, \dots, \numVertex,$
are suitable vertex matrices.

Second, based on this polytopic structure of the system, 
we consider polytopic LPV controllers of the form 
\begin{align}
\quad K_t^k &= \sum_{\indexVertex = 1}^{\numVertex} \coeffvertex \Kvertex,
\label{eq:K_t}
\end{align}
where $\Kvertex \in \real^{N \times N}, \indexVertex = 1, \dots, P,$ are vertex
matrices.  In Section \ref{sec:K_verteces}, we show how to design these vertex
matrices
so that the stabilizability and sparsity conditions are satisfied.

As we will show later, a polytopic LPV controller, with the ad-hoc sparsity
conditions imposed on each $\Kvertex$, for all $\indexVertex$, enables us to
divide the computation of $K_t^k, t \Tzmo$, for all $k$, into an offline step
and a distributed online one, thus making the optimal control algorithm fully
distributed.

\subsection{Polytopic LPV system and sparsity of the controller}
\label{sec:polytopic_representation}
First, we show how to compute vertex matrices and coefficients, respectively,
$\Avertex$ and $\coeffvertex$, for all $\indexVertex$, such that $A_t^k$ can be
written as in \eqref{eq:At_convex}-\eqref{eq:coefficients}.  The proposed
polytopic LPV representation of system \eqref{eq:system_z1} is a slightly
modified version of the one in \cite{millerioux2003polytopic}.
Let us consider the following assumption.
\begin{assumption}
\label{ass:A_t^k}
Every nonzero element of $A_t^k$ in \eqref{eq:system_z1} satisfies
$\at{(i,j)} \in [\amin{ij}, \amax{ij}]$, for all $ t \Tzmo$, with given
$\amin{ij} \in \real$ and $\amax{ij} \in \real$.  \oprocend
\end{assumption}

Then, let us associate an index $\indexTv = 1, \dots, \numTv$, where 
$\numTv := \sum_{i = 1}^N |\NN_i|$, to each nonzero element of $A_t^k$.
In the following, we will write $\at{(\indexTv)}$ when we use the index $s$
or $\at{(i,j)}$ when we use the index $i,j$. Moreover, we will use
$\amin{\indexTv}$ and  $\amax{\indexTv}$ to indicate the corresponding
$\amin{ij}$ and $\amax{ij}$, respectively.
\begin{proposition}
  Let $A_t^k$, $t \Tzmo$, satisfy Assumption~\ref{ass:A_t^k} for all $k$, then
  it can be written as in \eqref{eq:At_convex}-\eqref{eq:coefficients} with
  $\numVertex = 2\numTv$,
\begin{align}
\coeffvertex =
\begin{dcases}
\frac{1}{\numTv}\frac{\amax{\indexVertex}-\at{(\indexVertex)}}{\amax{\indexVertex}-\amin{\indexVertex}}
,\;  & \indexVertex = 1, \dots, \numTv,\\
\frac{1}{\numTv} \Big( 1 - \frac{\amax{\indexVertex-\numTv}-\at{(\indexVertex-\numTv)}}{\amax{\indexVertex-\numTv}-\amin{\indexVertex-\numTv}}\Big)
,\;  & \indexVertex = \numTv+1, \dots, 2\numTv,
\end{dcases}
\label{eq:coefficients_p}
\end{align}
and 
\begin{align}
\Avertex =
\begin{dcases}
\Amin_{\indexVertex}, \;  & \indexVertex = 1, \dots, \numTv,\\
\Amax_{\indexVertex-\numTv},\;  & \indexVertex = \numTv+1, \dots, 2\numTv,
\end{dcases}
\label{eq:A_p}
\end{align}
where $\Amin_{\indexTv} \in \real^{N \times N}$,
$\Amax_{\indexTv} \in \real^{N \times N}$, $\indexTv = 1, \dots, \numTv,$ are
matrices with all zeros except the elements indexed by $s$, which are
\begin{align*}
\AminElement_{\indexTv(\indexTv)}:= \numTv \;\! \amin{\indexTv} \quad \text{and}
  \quad \AmaxElement_{\indexTv(\indexTv)}:= \numTv \;\! \amax{\indexTv}. \tag*{$\square$}
\end{align*}
\end{proposition}

We now show how the sparsity condition on $K_t^k, \; t \Tzmo$, for all $k$, can be obtained by means of sparsity conditions imposed on each $\Kvertex$, for all $\indexVertex$.

\begin{proposition}
\label{prop:sparsity}
Let the vertex matrices $\Kvertex, \; \indexVertex = 1, \dots, \numVertex,$ satisfy, for all $k$, the sparsity condition
\end{proposition}
\begin{align}
\Kvertex \circ \sparseMatrixVertex^c = 0, \quad \indexVertex = 1, \dots, \numVertex,
\label{eq:sparsity_Kverteces}
\end{align} 
where $\sparseMatrixVertex^c := \boldsymbol{1} - \sparseMatrixVertex$ and
$\sparseMatrixVertex \in \real^{N \times N}, \; \indexVertex = 1, \dots, \numVertex$, is defined as a matrix, with the same sparsity of the corresponding $\Avertex$, such that 
\begin{align*}
\sparseMatrixVertexElement{(i,j)} =
\begin{cases}
0, \quad \text{if } \; \AvertexElement{(i,j)} = 0,\\ 
1, \quad \text{if } \; \AvertexElement{(i,j)} \neq 0.
\end{cases}
\end{align*}
Then, for all $k$, matrix $K_t^k, t \Tzmo$, as in \eqref{eq:K_t}, satisfies the sparsity condition in 
\eqref{eq:sparsity_distributed}. \oprocend

\subsection{Computation of vertex feedback matrices with sparsity and stability properties}
\label{sec:K_verteces}

In order to obtain stabilizing controller matrices,
we consider the necessary and sufficient condition 
for polytopic LPV systems, presented in
\cite{daafouz2001poly}.
This condition only guarantees the stabilizing property of the controller,
without taking into account the sparsity condition
\eqref{eq:sparsity_distributed}.

\begin{theorem}[\cite{daafouz2001poly}]
\label{eq:theorem_LMI}
System \eqref{eq:system_z1}, \eqref{eq:At_convex}, with output
$y_t \in \real^{N}$ such that $y_t = C \xlin_t + D \ulin_t, \; t \Tzmo$, where
$C \in \real^{N \times N}$ and $D \in \real^{N \times N}$ are given matrices, is
uniformly asymptotically stabilizable
with a given 
performance $\perfHinf$, 
by a control law \eqref{eq:control_law}, \eqref{eq:K_t}, if and only if there
exists matrices
$S_\indexVertex \in \real^{N \times N}, G_\indexVertex \in \real^{N \times N}$
and $R_\indexVertex \in \real^{N \times N}$,
$\indexVertex = 1, \dots, \numVertex$, satisfying for all
$\indexVertex = 1, \dots, \numVertex, \; q = 1, \dots, \numVertex$,
\begin{align}
\begin{split}
&
\small{\begin{bmatrix}
G_\indexVertex + G_\indexVertex^\top - S_\indexVertex
 & 0^\top &\!\!\! (\tilde{A}_\indexVertex G_\indexVertex - B R_\indexVertex)^\top &\!\!\!\!\! (C G_\indexVertex - D R_\indexVertex)^\top\\
0 & \perfHinf I & 0 &0\\
\Avertex G_\indexVertex - B R_\indexVertex & 0 & S_q & 0\\
C G_\indexVertex - D R_\indexVertex & 0 & 0 & \perfHinf I
\end{bmatrix}} > 0,\\[0.5ex]
& \quad S_\indexVertex = S_\indexVertex^\top, \; S_\indexVertex > 0.
\label{eq:LMIs}
\end{split}
\end{align}
\normalsize
The stabilizing controller $K_t^k, t \Tzmo,$ is given by \eqref{eq:K_t} with 
$\Kvertex = R_\indexVertex G_\indexVertex^{-1}$, for all $p$.
\oprocend
\end{theorem}

To obtain both sparse and stabilizing controllers $K_t^k$, we thus need to
design vertex controllers $\Kvertex$ with these properties,
i.e., $\Kvertex \circ \sparseMatrixVertex^c = 0$ and
$\Kvertex = R_\indexVertex G_\indexVertex^{-1}$, for all $\indexVertex$.
Notice that,
\begin{align}
R_\indexVertex G_\indexVertex^{-1} \circ \sparseMatrixVertex^c = 0,
\label{eq:nonconvex_condition}
\end{align}
is a nonconvex constraint, so that the controller computation turns out to be a
challenging problem.

Thus, by taking inspiration from \cite{fardad2014design}, we propose a novel
algorithm (Algorithm \ref{alg:algorithm_Kvertices}) for the computation of the
vertex matrices $\Kvertex, \; \indexVertex = 1, \dots, \numVertex$.
We consider the conditions \eqref{eq:LMIs} together with the sparsity condition
\eqref{eq:nonconvex_condition}, where the matrices $G_\indexVertex$, for all
$\indexVertex$, are replaced by an estimate $\hat{G}_\indexVertex$, thus getting
the conditions \eqref{eq:problem_convex}.  We initialize
$\hat{G}_\indexVertex = I, \; \indexVertex = 1, \dots, \numVertex$, where $I$ is
the identity matrix. Then, at iteration $h$, we compute the matrices
$G_\indexVertex^h, R_\indexVertex^h, S_\indexVertex^h,$ for all $\indexVertex$,
satisfying \eqref{eq:problem_convex} and we set
$\hat{G}_\indexVertex = G_\indexVertex^h$, until we reach a good approximation
of the matrices $G_\indexVertex$ via $\hat{G}_\indexVertex$, for all
$\indexVertex$, and the sparsity condition is satisfied with a given
tolerance. When the latter conditions are reached, we get the vertex matrices of
the controller as $\tilde{K}_p = R_p^h G_p^{h^{-1}}, \; p = 1, \dots, P$.
By using this strategy, we get a sequence of convex feasibility problems that
can be solved at each iteration via available numerical toolboxes.

\begin{algorithm}[H]
\caption{Computation of $\Kvertex, \; \indexVertex = 1, \dots, \numVertex,$}
\label{alg:algorithm_Kvertices}
\begin{algorithmic}
\small
\Require 
$\Avertex, \hat{G}_\indexVertex = I$, for all $\indexVertex$, $B, C, D, \perfHinf, \epsilon > 0$
\vspace{0.1cm}

\For{$h = 0, 1, \ldots$}\vspace{0.2cm}
\State find $G_\indexVertex^h, R_\indexVertex^h, S_\indexVertex^h, \; \indexVertex = 1, \dots, \numVertex,$ satisfying
\begin{align}
\begin{split}
& \!\!\!
\begin{bmatrix}
G_\indexVertex + G_\indexVertex^\top -S_\indexVertex & 0^\top &\!\!\! (\Avertex G_\indexVertex - B R_\indexVertex)^\top \!\!\!&\!\!\! (C G_\indexVertex - D R_\indexVertex)^\top\\
0 & \perfHinf I & 0 &0\\
\Avertex G_\indexVertex - B R_\indexVertex & 0 & S_q & 0\\
C G_\indexVertex - D R_\indexVertex & 0 & 0 & \perfHinf I
\end{bmatrix} > 0,\\
& S_\indexVertex = S_\indexVertex^\top, \; S_\indexVertex > 0,\\
& R_\indexVertex \hat{G}_\indexVertex^{-1} \circ \sparseMatrixVertex^c = 0, \quad \indexVertex = 1, \dots, \numVertex, \quad q = 1, \dots, \numVertex
\label{eq:problem_convex}
\end{split}
\end{align}
\If {$||\hat{G}_\indexVertex - G_\indexVertex^h|| \! < \epsilon, ||R_\indexVertex^h G_\indexVertex^{h^{-1}} \!\!\!\circ \sparseMatrixVertex^c|| < \epsilon, \; \text{for all } \indexVertex,$}\vspace{0.05cm}
	\State $\Kvertex = R_\indexVertex^h G_\indexVertex^{h^{-1}}\!\!\!\!, \; \text{for all } \indexVertex,$ and \Return 
\Else
	\State $\text{ set }\hat{G}_\indexVertex = G_\indexVertex^h, \; \text{for all } \indexVertex$	
\EndIf
\EndFor
\end{algorithmic}
\end{algorithm}
\normalsize

\section{Fully-Distributed Algorithm\\ for Nonlinear Optimal Control}
\label{sec:fully-distributed}
In this section we show how to obtain a fully-distributed version of the
cloud-assisted distributed algorithm, \cite{spedicato2018necsys}, recalled in
Section \ref{sec:cloud-assisted-algorithm}.

Once vertex matrices $\Kvertex$, for all $p$, are computed by means
of Algorithm \ref{alg:algorithm_Kvertices}, a sparse stabilizing
$K_t^k, \; t \Tzmo$, can be obtained by means of \eqref{eq:K_t}.
We now show how agent $i$ can compute the elements $k_{t(j,i)}^k, k_{t(i,j)}^k, \; j \in \NN_i$, via the only information of neighbors. 
By defining 
 \begin{align}
 \Kmin_\indexTv := \Kv_\indexTv,\quad
 \Kmax_\indexTv := \Kv_{\indexTv+\numTv}, \quad \indexTv = 1,\dots,\numTv, 
 \label{eq:defK}
 \end{align}
and using \eqref{eq:K_t} and \eqref{eq:coefficients_p}, the controller $K_t^k$ can be written as 
 \begin{align*}
 K_t^k \!&= \!\!\sum_{\indexTv = 1}^\numTv\! \Big(\frac{1}{\numTv} 
 \Big(\frac{\amax{\indexTv} -\at{(\indexTv)}}{\amax{\indexTv} - \amin{\indexTv}} \Big) \Kmin_{\indexTv} \!+\!
 \frac{1}{\numTv} \Big( \!1 - \!
 \frac{\amax{\indexTv} - \at{(\indexTv)}}{\amax{\indexTv} - \amin{\indexTv}}
 \Big) \! \Kmax_{\indexTv} \! \Big).
 \end{align*}
Comparing the definitions of $\Amin_{\indexTv}$ (resp. $\Amax_{\indexTv}$)
with $\Kmin_{\indexTv}$ (resp. $\Kmax_{\indexTv}$), it follows that they have
the same sparsity.
In particular, each one of them has only one nonzero element, specifically the
element $\KminElement_{\indexTv(\indexTv)}$
(resp. $\KmaxElement_{\indexTv(\indexTv)}$).
The nonzero elements $\KminElement_{s(s)}, \KmaxElement_{s(s)}$ for all
$\indexTv = 1, \dots, \numTv$ can be equivalently indexed by $i,j$ instead of
the $s$ and written as $\KminElement_{ij(i,j)}, \KmaxElement_{ij(i,j)}$,
correspondingly. They can be computed as
	\begin{align*}
	\small
	k_{t(i,j)}^{k} &\!\!=\!\! 
	\frac{1}{S}\! \Bigg(\!\!\!
	\bigg(\!\!\frac{\amax{ij} -\at{(i,j)}}{\amax{ij} - \amin{ij}} \!\!\bigg) \KminElement_{ij(i,j)} \!\!+ \!\!
	\bigg( \!1\!\! - \!
	\frac{\amax{ij} - \at{(i,j)}}{\amax{ij} - \amin{ij}}\!\!
	\bigg)\! \KmaxElement_{ij(i,j)} \!\!\Bigg),
	\end{align*}
	for all $i = 1, \dots, N$, $j \in \NN_i$.
 Clearly, each agent $i$ can compute this expression in a completely distributed
 way.

In the next table (Algorithm \ref{alg:distributed_algorithm}) our fully-distributed algorithm is presented from the perspective of agent $i$.
Differently from the cloud-assisted distributed
algorithm, agent $i$ computes by means of
\eqref{eq:K_t_ji_alg} the gains $k_{t(j,i)}^k, \; j \in \NN_i$, that are used for
the descent direction in \eqref{eq:descent_zp}, and by means of
\eqref{eq:K_t_ij_alg} the gains $k_{t(i,j)}^k, \; j \in \NN_i,$ used for the
trajectory update in \eqref{eq:closed_loop_update}. No central unit is used for the controller computation. 

\begin{algorithm}[H]
\caption{Fully-distributed version of Algorithm \ref{alg:cloud_assisted}}
\label{alg:distributed_algorithm}
\begin{algorithmic}
\Require $x_{j,t}^{0}$, $u_{i,t}^{0}$, for all $t$, $j \in \NN_i$, with $(x^{0}
u^{0})$ a trajectory of \eqref{eq:dynamics}, constant step-size $\beta$, $\amin{ij}, \amax{ij}$, $\KminElement_{ij(i,j)}, \KmaxElement_{ij(i,j)},
\amin{ji}, \amax{ji}, \KminElement_{ji(j,i)}, \KmaxElement_{ji(j,i)}$, $j \in \NN_i$
\vspace{0.1cm}

\For{$k = 0, 1, 2 \ldots$}
\State set $p_{i,\finalTime}^k = \nabla m_i(x_{i,\finalTime}^k)$
\For{$t = \finalTime-1, \dots, 0$}
\State compute $\udir_{i,t}^k$ via \eqref{eq:descent_v}
\State receive $\!\at{(j,i)}, b_{(j,j)}, \udir_{j,t}^k, \ell_{u,j,t}^k, p_{j, t+1}^k$, $j\!\in\!\NN_i \!\setminus \!\{i\}$
\State
\vspace{-0.2cm}
\State compute, for all $j \in \NN_i$,
\small{
\begin{align}
k_{t(j,i)}^{k} &\!\!=\!\!  
\frac{1}{S} \!\Bigg(\!\!
\!\bigg(\!\!\frac{\amax{ji} -\at{(j,i)}}{\amax{ji} - \amin{ji}} \!\!\bigg) \KminElement_{ji(j,i)} \!\!+ \!\!
\bigg( \!1 \!- \!
\frac{\amax{ji} - \at{(j,i)}}{\amax{ji} - \amin{ji}}
\!\!\bigg) \KmaxElement_{ji(j,i)} \!\!\Bigg),
\label{eq:K_t_ji_alg}
\end{align}
}
\State compute $\xdir_{i,t}^k, p_{i,t}^k$ via \eqref{eq:descent_zp}
\vspace{0.1cm}
\EndFor
\vspace{-0.2cm}
\For{$t = 0, 1, \ldots, \finalTime-1$} 
\State receive  $\xdir_{j,t}^k$, $x_{j,t}^{k+1}$,  $j \in \NN_i \setminus \{i\}$
\State \multiline{compute $\alpha_{j,t}^{k+1}, \; j \in \NN_i, \; \mu_{i,t}^{k+1}\!$ via \eqref{eq:curve_update}, with $\beta^k = \beta$}
\State compute, for all $j \in \NN_i$,
\small{
\begin{align}
k_{t(i,j)}^{k} &\!\!=\!\! 
\frac{1}{S}\! \Bigg(\!\!\!
\bigg(\!\!\frac{\amax{ij} -\at{(i,j)}}{\amax{ij} - \amin{ij}} \!\!\bigg) \KminElement_{ij(i,j)} \!\!+ \!\!
\bigg( \!1\!\! - \!
\frac{\amax{ij} - \at{(i,j)}}{\amax{ij} - \amin{ij}}\!\!
\bigg) \KmaxElement_{ij(i,j)} \!\!\Bigg),
\label{eq:K_t_ij_alg}
\end{align}
}
\State compute $u_{i,t}^{k+1}, x_{i,t+1}^{k+1}$ via \eqref{eq:closed_loop_update}
\EndFor
\EndFor
\end{algorithmic}
\end{algorithm}

\section{Simulations}
\label{sec:numerical_computations}
In this section we present a numerical example to show the effectiveness of the
proposed strategy for the design of sparse stabilizing controllers.
We consider a multi-agent system implementing a (distributed) formation control
law based on virtual potential functions, see, e.g.,
\cite{rozenheck2015proportional}, and equip agents with an additional input.
Specifically, the local state function of the nonlinear dynamics over
graph \eqref{eq:dynamics}, for each agent $i=1,\ldots, N$, is given by
\begin{align*}
 &f_i(x_{\NN_{i},t},u_{i,t}) = x_{i,t}\nonumber \\
 & \quad - T_s(\| x_{i,t} -x_{i+1,t}\|^2 - d_{i,i+1}^2)(x_{i,t}-x_{i+1,t}) \nonumber \\ 
 & \hspace{4mm} -T_s (\| x_{i,t} -x_{i-1,t}\|^2 - d_{i,i-1}^2)(x_{i,t} -x_{i-1,t}) + T_s c\,u_{i,t},
\end{align*}
where $x_{i,t} \in \real^2$ is the position vector and $u_{i,t} \in \real^2$ is
the input vector of agent $i$ at time instant $t$. The parameters $d_{i,i+1}$
and $d_{i,i-1}$ represent distances between agents $i$ and $i+1$, and agents $i$
and $i-1$, respectively, in the desired formation. The parameter $c$ is an input
coefficient, while $T_s$ is the sampling time used for the discretization in
time. Here we are considering $N=6$ agents interacting according to a cycle graph, so
that $x_{0,t} = x_{N,t}$ and $x_{N+1,t} = x_{1,t}$.
We set $T_s=10^{-2}$ and
$c=10$. The target formation of the multi-agent system is a hexagon with side
length of $4$m, so that we set $d_{i,i+1} = d_{i,i-1} = 4$m.
We design the sparse polytopic LPV controller by using $\nu = 0.05$, $C=I$ and
$D=10^{-5} I$.

We generated the vertex matrices according to
Algorithm~\ref{alg:algorithm_Kvertices}, chose a trajectory to be stabilized (by
simply integrating the open-loop system), and thus generated the (time-varying)
stabilizing feedback $K_t$ for the given trajectory.

\begin{figure}[h]
\hspace{2cm}
	\includegraphics[scale = 0.3]{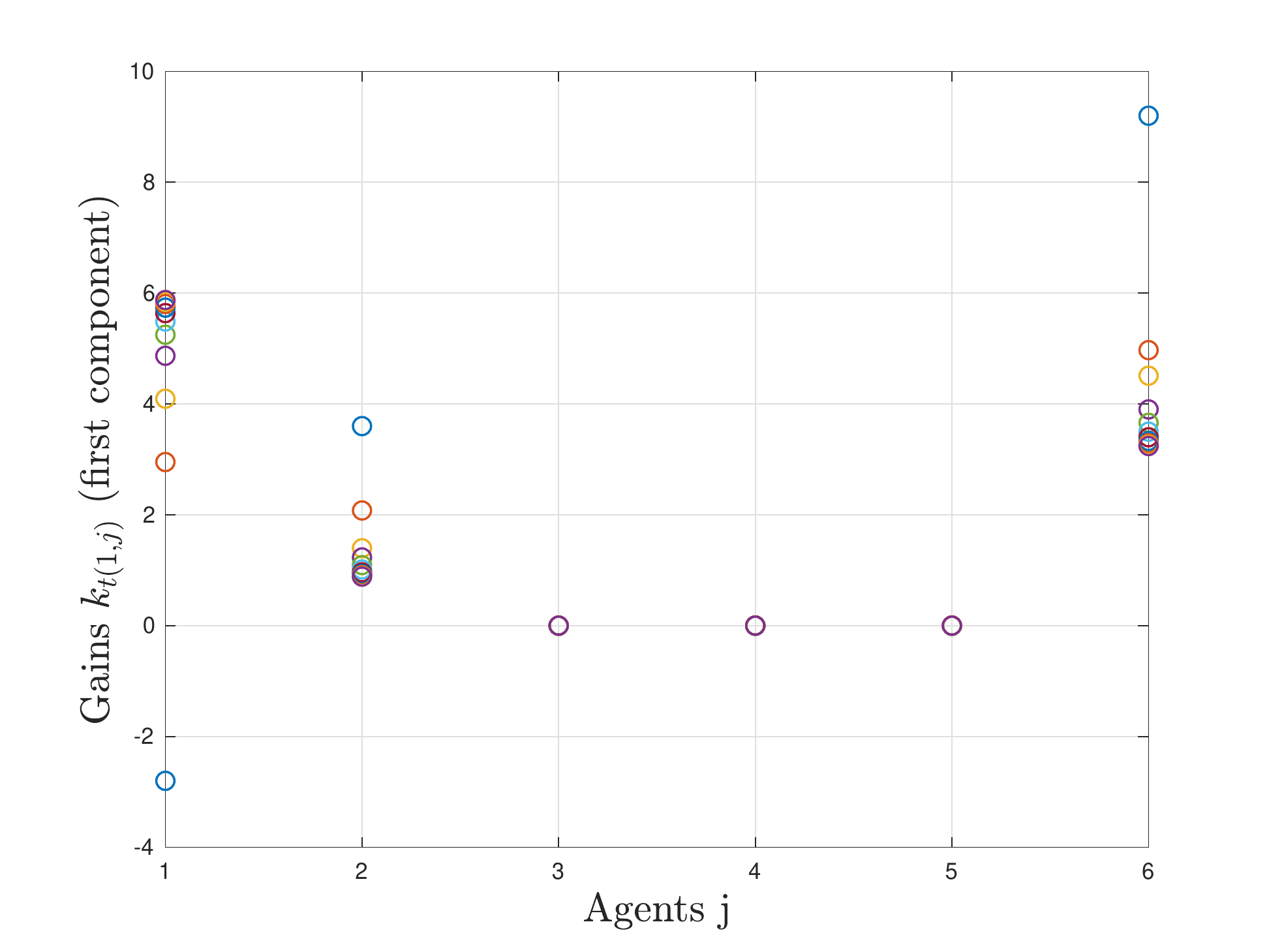}
	\includegraphics[scale = 0.3]{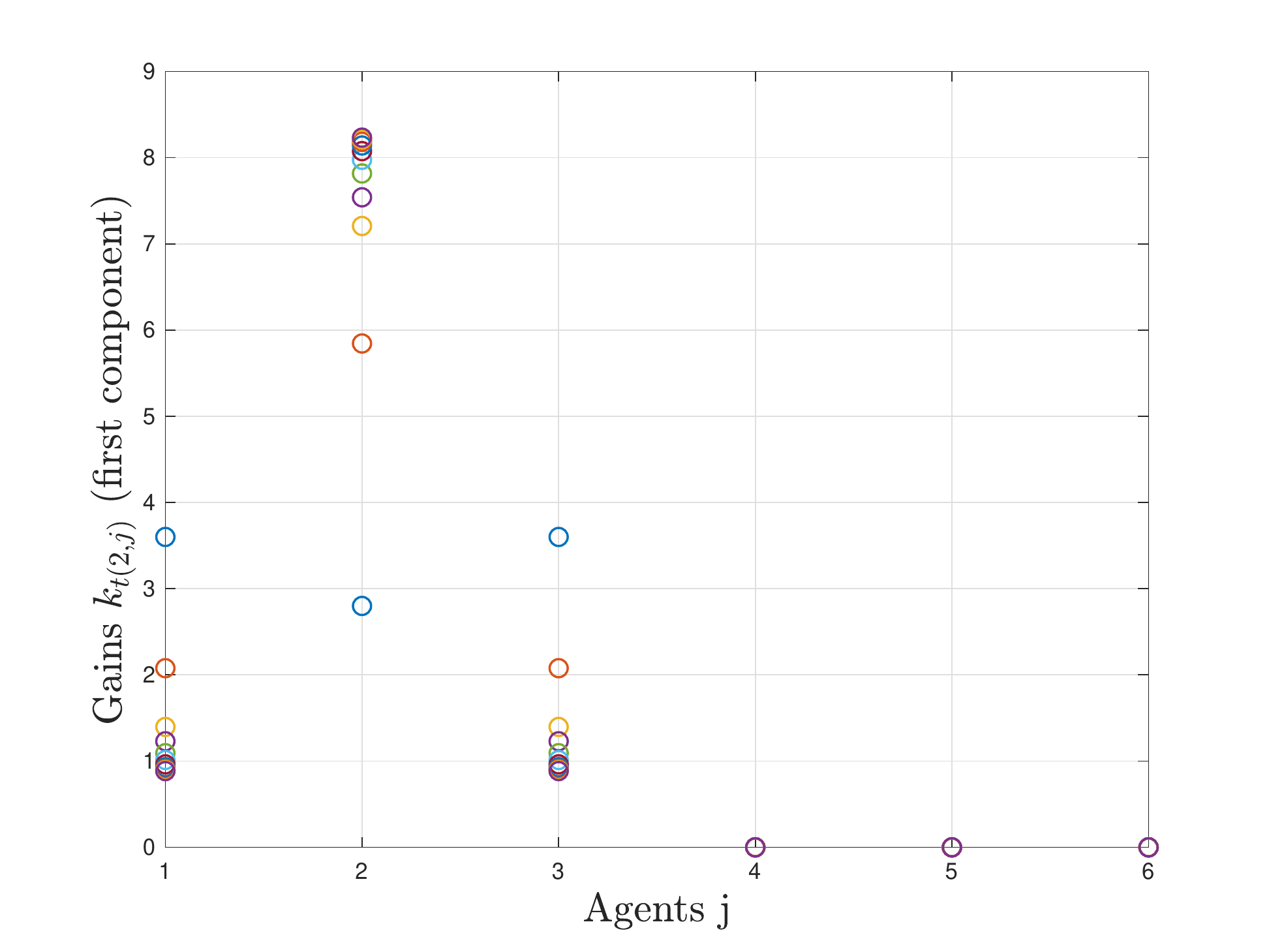}
	\caption{Gain values $k_{t(i,j)}$ (first component) for agents $i=1$
          (left) and $i=2$ (right) plotted with respect to $j=1,\ldots,6$.}
	\label{fig:gainplot}
\end{figure}

In Figure~\ref{fig:gainplot}, we depicted the values of $K_t$ for the first (position)
component of agents $1$ and $2$. In particular, to highlight the sparsity of the
controller, for each agent $i=1,2$ we plotted the values $k_{t(i,j)}$ (first
component), $t=1,\ldots,10$, (circles of different colors) with respect to
$j = 1,\ldots,6$. As expected the weights $k_{t(i,j)}$ with $j\notin \NN_i$ are
zero for all $t$.

\begin{figure}[h]
	\centering
	\includegraphics[scale = 0.5]{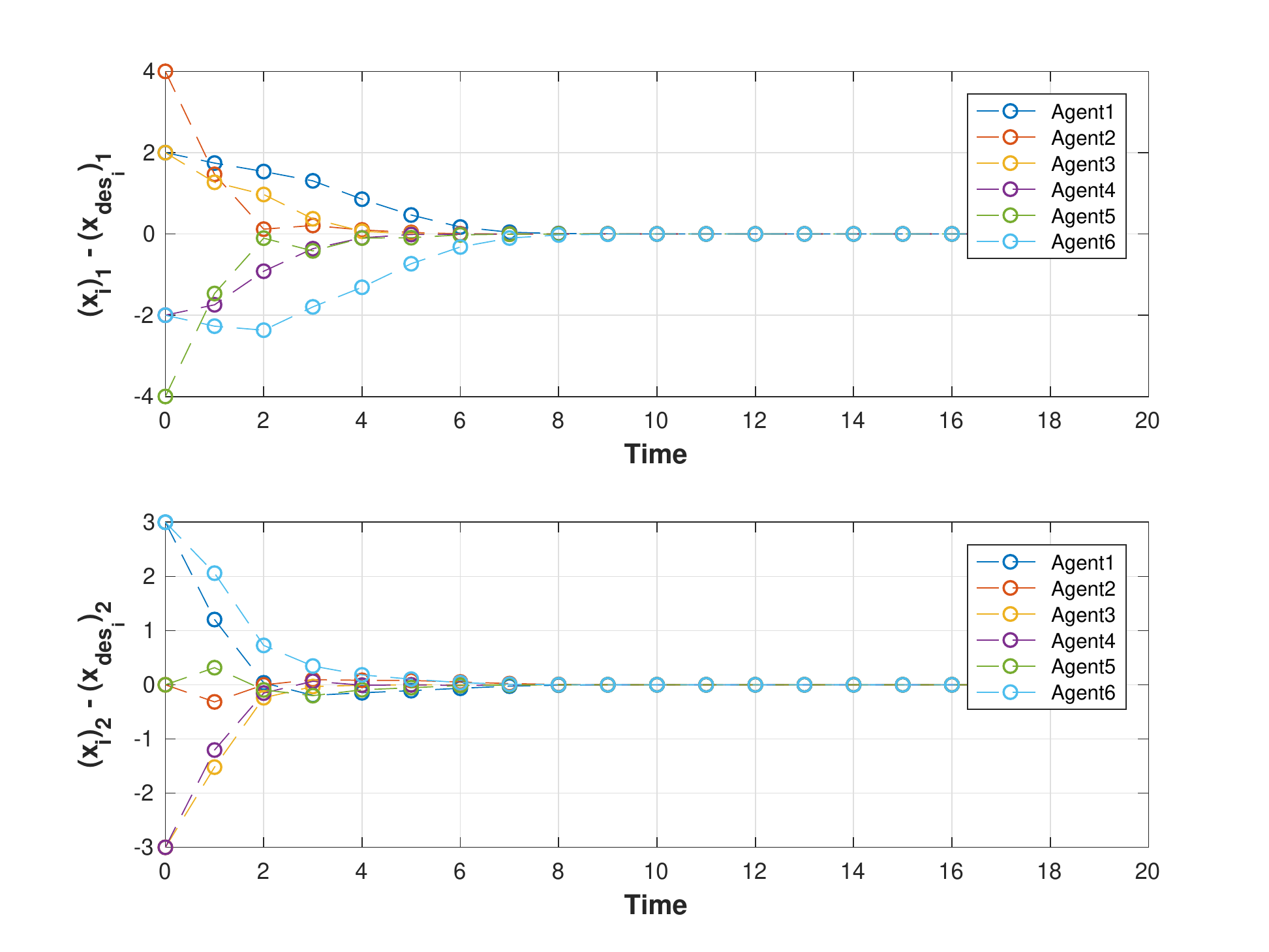}
	\caption{Error with respect to the desired trajectory $x_{\text{des}_i}$
          for all agents $i=1,\ldots,6$.}
	\label{fig:errorplot}
\end{figure}

In Figure \ref{fig:errorplot}, we show the evolution of the position error with
respect to the desired trajectory $x_{\text{des}_i}$ for the two components
$(x_i)_1$ and $(x_i)_2$ of all agents $i = 1,\ldots, N$.

\section{Conclusions}
In this paper we have proposed a fully-distributed strategy to solve nonlinear
optimal control problems over networks. The strategy extends the one proposed in
our previous work~\cite{spedicato2018necsys} in which a cloud was used together
with distributed computation. The main distinctive feature of the new algorithm
is the design of a sparse stabilizing controller allowing agents to stabilize
system trajectories (at each iteration of the optimal control algorithm) in a
distributed way. By relying on polytopic LPV systems, we proposed a two step
procedure for the controller design. First, we proposed an iterative algorithm,
to be performed offline, to compute ``stabilizing'' vertex feedback matrices
satisfying nonconvex sparsity constraints. Second, we showed how at each
iteration nodes can compute feedback gains stabilizing the current trajectory in
a fully-distributed way. The controller was tested in simulation on a
multi-robot formation control problem.

\bibliographystyle{IEEEtran}
\bibliography{bibliography}

\end{document}